\ifx\shlhetal\undefinedcontrolsequence\let\shlhetal\relax\fi

\documentclass[11pt]{amsart}

\usepackage{amsmath}
\usepackage{amssymb}

\newtheorem{theorem}{Theorem}[section]
\newtheorem{claim}[theorem]{Claim}

\newtheorem{proposition}[theorem]{Proposition}
\newtheorem{corollary}[theorem]{Corollary}

\theoremstyle{definition}
\newtheorem{definition}[theorem]{Definition}

\theoremstyle{remark}
\newtheorem{remark}[theorem]{Remark}

\newcount\skewfactor
\def\mathunderaccent#1#2 {\let\theaccent#1\skewfactor#2
\mathpalette\putaccentunder}
\def\putaccentunder#1#2{\oalign{$#1#2$\crcr\hidewidth
\vbox to.2ex{\hbox{$#1\skew\skewfactor\theaccent{}$}\vss}\hidewidth}}
\def\name{\mathunderaccent\tilde-3 }

\def\smallbox#1{\leavevmode\thinspace\hbox{\vrule\vtop{\vbox
   {\hrule\kern1pt\hbox{\vphantom{\tt/}\thinspace{\tt#1}\thinspace}}
   \kern1pt\hrule}\vrule}\thinspace}

\newcommand{\cf}{{\rm cf}}

\newcommand{\Dom}{{\rm Dom}}

\newcommand{\then}{{\underline{then}}}
\newcommand{\Then}{{\underline{Then}}}

\def\qedref#1{$\qed_{\ref{#1}}$}

\title{Combinatorial aspects of the splitting number}
\author{Shimon Garti}
\address{Institute of Mathematics
 The Hebrew University of Jerusalem
 Jerusalem 91904, Israel}
\email{shimon.garty@mail.huji.ac.il}

\author{Saharon Shelah}
\address{Institute of Mathematics
 The Hebrew University of Jerusalem
 Jerusalem 91904, Israel
 and  Department of Mathematics
 Rutgers University
 New Brunswick, NJ 08854, USA}
\email{shelah@math.huji.ac.il}
\urladdr{http://www.math.rutgers.edu/\char`\~shelah}

\thanks{First typed: December 2010 \\
Research supported by the United States-Israel Binational
Science Foundation. Publication 962 of the second author}

\begin{document}

\begin{abstract}
This paper deals with the splitting number $\mathfrak{s}$ and polarized partition relations. In the first section we define the notion of strong splitting families, and prove that its existence is equivalent to the failure of the polarized relation $\binom{\mathfrak{s}}{\omega} \rightarrow \binom{\mathfrak{s}}{\omega}^{1,1}_2$. We show that the existence of a strong splitting family is consistent with ZFC, and that the strong splitting number equals the splitting number, when it exists. Consequently, we can put some restriction on the possibility that $\mathfrak{s}$ is singular. In the second section we deal with the polarized relation under the weak diamond, and we prove that the strong polarized relation $\binom{2^\omega}{\omega} \rightarrow \binom{2^\omega}{\omega}^{1,1}_2$ is consistent with ZFC, even when $\cf(2^\omega)=\aleph_1$ (hence the weak diamond holds).
\end{abstract}

\maketitle

\newpage

\section{introduction}
This paper deals with two problems. The first is a topological one, and it deals with strong splitting families in $\mathcal{P}(\omega)$. the second is a combinatorial one, and related to the polarized partition relation.

Let us start with the topological problem. A family $\bar{S} = \{ S_\alpha : \alpha < \kappa \} \subseteq \mathcal{P}(\omega)$ is splitting if for every $B \in [\omega]^\omega$ there exists an ordinal $\alpha < \kappa$ so that $|B \cap S_\alpha| = |B \setminus S_\alpha| = \aleph_0$. In this case we say that $S_\alpha$ splits $B$. The cardinal invariant $\mathfrak{s}$, the splitting number, is defined as the minimal cardinality of a splitting family. A good source for information about $\mathfrak{s}$ (as well as other basic cardinal invariants on the continuum) is van Dowen (in \cite{MR776622}). We follow, in this paper, his terminology.

Notice that the existence of one ordinal $\alpha$ such that $S_\alpha$ splits $B$ is enough for this definition. We may ask, further, if one can find a family of subsets of $\omega$ so that for each $B \in [\omega]^\omega$ almost every set among the sets in the splitting family splits $B$. This is the property of strong splitting families, and our topological problem is whether such a property is possible.

For splitting families in the common sense, one can always take the collection of all the subsets of $\omega$. But this does not work for strong splitting families. On one hand, we need enough sets (in the family that we try to create) so that every $B \in [\omega]^\omega$ is split. On the other hand, we must be careful not to take too many sets, since otherwise we will have some $B \in [\omega]^\omega$ which is included in a lot of sets.

Let us describe the combinatorial problem. The balanced polarized relation $\binom{\alpha}{\beta} \rightarrow \binom{\gamma}{\delta}^{1,1}_2$ asserts that for every coloring $c : \alpha \times \beta \rightarrow 2$ there are $A \subseteq \alpha$ and $B \subseteq \beta$ so that ${\rm otp} (A) = \gamma, {\rm otp} (B) = \delta$ and $c \upharpoonright (A \times B)$ is constant. This relation was first introduced in papers of Erd\"os, Hajnal and Rado (see \cite{MR0202613}, and \cite{MR0081864}). A good reference for the basic facts about this relation is \cite{williams}. 

If $\alpha = \gamma$ and $\beta = \delta$, we name this relation as a strong polarized relation. Our question is whether the strong relation $\binom{\mathfrak{s}}{\omega} \rightarrow \binom{\mathfrak{s}}{\omega}^{1,1}_2$ holds. As we shall see, the topological problem above is deeply connected to this combinatorial question. In fact, our ability to solve the combinatorial part enables us to give an answer to the topological problem.

In the second section we deal with the polarized relation on the continuum. The starting point is the negative result of Erd\"os and Rado, that $\binom{\omega_1}{\omega} \nrightarrow \binom{\omega_1}{\omega}_2^{1,1}$ under the continuum hypothesis. Of course, the positive relation is also consistent (e.g., under the PFA). 

We have asked whether the correct generalization of the negative result is $\binom{2^\omega}{\omega} \nrightarrow \binom{2^\omega}{\omega}^{1,1}_2$. It was proved in \cite{gash964} that the positive relation $\binom{2^\omega}{\omega} \rightarrow \binom{2^\omega}{\omega}^{1,1}_2$ is consistent with ZFC. But in that paper, $\cf(2^\omega) \geq \aleph_2$, and one of the referees suggested to consider the possibility that $\cf(2^\omega)=\aleph_1 \Rightarrow \binom{2^\omega}{\omega} \nrightarrow \binom{2^\omega}{\omega}^{1,1}_2$. We shall prove the converse.

Let us try to explain the background of this interesting suggestion of the referee. If $\Diamond_{\aleph_1}$ holds, then $\binom{2^\omega}{\omega} \nrightarrow \binom{2^\omega}{\omega}^{1,1}_2$. If $\cf(2^\omega)=\aleph_1$ then we know that a weak version of the diamond holds. This version is the so-called weak diamond, and we denote it by $\Phi_{\aleph_1}$. It follows that this principle is equivalent to the cardinal assumption $2^{\aleph_0} < 2^{\aleph_1}$, and since $\cf(2^\omega)=\aleph_1 \Rightarrow 2^{\aleph_0} < 2^{\aleph_1}$ one might guess that under this assumption we shall get $\binom{2^\omega}{\omega} \nrightarrow \binom{2^\omega}{\omega}^{1,1}_2$. As we shall see, the weak diamond does not imply this negative result (in contrary to the full diamond).
We also show that $\Phi_{\aleph_1}$ is consistent with the positive relation $\binom{\omega_1}{\omega} \rightarrow \binom{\omega_1}{\omega}_2^{1,1}$.

We try to use standard notation. The combinatorial notation is due to \cite{MR795592}. We save the letter $H$ for monochromatic sets, when possible. The symbol $A \subseteq^* B$ means that $|A \setminus B| < \aleph_0$. We use $\kappa, \lambda, \mu, \tau$ as cardinals, and $\alpha, \beta, \gamma, \delta, \varepsilon, \zeta$ as ordinals. $n$ is a finite ordinal, and $\omega$ is the first infinite ordinal. We use $\mathfrak{c}$ and $2^{\aleph_0}$ interchangeably. The second section employs forcing arguments. 
Some benighted people reverse the natural order in forcing relations.
We indicate that $p \leq q$ means (in this paper) that $q$ gives more information than $p$ in forcing notions. For background in forcing (including the notation we adhere to) we suggest \cite{MR1623206}. For background in pcf theory (e.g., the covering numbers which appear at the end of the first section) the reader may consult the monograph \cite{MR1318912}.

\newpage

\section{Splitting families and the polarized relation}
\par \noindent Let us define the strong splitting property:

\begin{definition}
\label{ssf}
Strong splitting. \newline 
Let $\mathcal{F} = \{ S_\alpha : \alpha < \kappa \}$ be a family of subsets of $\omega$, and assume $\mathcal{F}$ is a splitting family. \newline 
For $B \in [\omega]^\omega$ set $\mathcal{F}_B = \{ S_\alpha : (B \subseteq^* S_\alpha)$ or $(B \subseteq^* \omega \setminus S_\alpha) \}$. \newline 
$\mathcal{F}$ is a strong splitting family if $|\mathcal{F}_B| < \mathfrak{s}$ for every $B \in [\omega]^\omega$. \newline 
$\mathcal{F}$ is a very strong splitting family if $|\mathcal{F}_B| < |\mathcal{F}|$ for every $B \in [\omega]^\omega$.
\end{definition}

\begin{remark}
\label{ttrrivial}
If $\mathcal{F} = \{ S_\alpha : \alpha < \kappa \}$ is a strong splitting family, then $\kappa \geq \mathfrak{s}$, since a strong splitting family is, in particular, a splitting family.
\end{remark}

We start with a claim that draws a connection of double implication between the topological question of strong splitting families, and the strong polarized relation with respect to the splitting number: 

\begin{claim}
\label{eequiv}
The equivalence claim.
\begin{enumerate}
\item [$(a)$] A strong splitting family (in $\mathcal{P}(\omega)$) exists iff $\binom{\mathfrak{s}}{\omega} \nrightarrow \binom{\mathfrak{s}}{\omega}^{1,1}_2$.
\item [$(b)$] A very strong splitting family of cardinality $\mu$ exists iff $\binom{\mu}{\omega} \nrightarrow \binom{\mu}{\omega}^{1,1}_2$.
\end{enumerate}
\end{claim}

\par \noindent \emph{Proof}. \newline 
We prove part $(a)$, and the same argument gives also part $(b)$.
Suppose $\binom{\mathfrak{s}}{\omega} \rightarrow \binom{\mathfrak{s}}{\omega}^{1,1}_2$, and assume toward contradiction that $\mathcal{F} = \{ S_\alpha : \alpha < \kappa \}$ is a strong splitting family.
We define a coloring $c_{\mathcal{F}} \equiv c : \mathfrak{s} \times \omega \rightarrow 2$ as follows:

$$
c(\alpha, n) = 0 \Leftrightarrow n \in S_\alpha
$$

This is done for every $\alpha < \mathfrak{s}$ and every $n \in \omega$. Since $\binom{\mathfrak{s}}{\omega} \rightarrow \binom{\mathfrak{s}}{\omega}^{1,1}_2$, there are $H_0 \in [\mathfrak{s}]^{\mathfrak{s}}$ and $H_1 \in [\omega]^\omega$ so that $c \upharpoonright (H_0 \times H_1)$ is constant. Without loss of generality, $\alpha \in H_0, n \in H_1 \Rightarrow c(\alpha, n) = 0$.

On one hand, $\mathcal{F}$ is strong splitting, so $|\mathcal{F}_{H_1}| < \mathfrak{s}$. On the other hand, if $\alpha \in H_0$ then $H_1 \subseteq S_\alpha$ (since $c(\alpha, n) = 0$ for every $n \in H_1$, and by the definition of $c$). So clearly, $H_1 \subseteq^* S_\alpha$ for every $\alpha \in H_0$. But $|H_0| = \mathfrak{s}$, so $|\mathcal{F}_{H_1}| \geq \mathfrak{s}$, a contradiction.

The opposite implication is similar. Suppose there is no strong splitting family, aiming to show that $\binom{\mathfrak{s}}{\omega} \rightarrow \binom{\mathfrak{s}}{\omega}^{1,1}_2$ holds. Given a coloring $c : \mathfrak{s} \times \omega \rightarrow 2$, we wish to find a monochromatic cartesian product with the desired cardinalities. For $\ell = 0,1$ and for every $\alpha < \mathfrak{s}$, set $S^\ell_\alpha = \{ n \in \omega : c(\alpha, n) = \ell \}$. Now, let $\mathcal{F}_c \equiv \mathcal{F}$ be the following family:

$$
\{ S^\ell_\alpha : \alpha < \mathfrak{s}, \ell = 0,1 \}
$$

By our assumption (in this direction), $\mathcal{F}$ is not a strong splitting family. Choose a witness, i.e., a set $B \in [\omega]^\omega$ such that $|\mathcal{F}_B| \geq \mathfrak{s}$. Collect the ordinals of $\mathcal{F}_B$, i.e., let $H_B$ be the set $\{ \alpha < \mathfrak{s} : \exists \ell \in \{0,1\}, S^\ell_\alpha \in \mathcal{F}_B \}$. Since $|\mathcal{F}_B| \geq \mathfrak{s}$ and $\ell$ ranges just over two values, we may assume (without loss of generality) that $\ell = 0$.

By a similar argument, we may assume that $B \subseteq^* S^0_\alpha$ (and not in its complement) for every $\alpha < \mathfrak{s}$. Moreover, we can replace (again, without loss of generality) the relation $\subseteq^*$ by $\subseteq$. This is justified by the fact that $\cf(\mathfrak{s}) > \aleph_0$ upon noticing that $|\omega^{< \omega}| = \aleph_0$. Consequently, one can choose a finite set of natural numbers which occurs $\mathfrak{s}$-many times as the discrepancy between $\subseteq^*$ and $\subseteq$. Now, remove this set from $B$, and we still have an infinite monochromatic set as required.
The last step is to observe that $c \upharpoonright (H_B \times B) \equiv 0$. Since the coloring $c$ was arbitrary, we are done.

\hfill \qedref{eequiv}

\begin{remark}
\label{ccohen}
If we add $\lambda=\lambda^{\aleph_0}$ Cohen reals, so $\mathfrak{c}=\lambda$, then for every $\mu\in(\aleph_0,\lambda]$ we have $\binom{\mu}{\omega} \nrightarrow \binom{\mu}{\omega}^{1,1}_2$. Moreover, $\binom{\mu}{\omega} \nrightarrow \binom{\aleph_1}{\omega}^{1,1}_2$.
\end{remark}

\par \noindent \emph{Proof}. \newline 
Let $\mathbb{Q}$ be the forcing notion which adds $\lambda$-many Cohen reals,  $\langle g_\alpha:\alpha<\lambda\rangle$. Define $\name{c}(\alpha,n)=\name{g}_\alpha(n)$. We claim that this coloring demonstrates the negative relation to be proved.

Towards contradiction assume that $\name{A}\in[\mathbb{N}]^{\aleph_0}, \name{B}\in[\mu]^{\aleph_1}$, and there exists a condition $p_0$ so that $p_0\Vdash_{\mathbb{Q}} \name{c}\upharpoonright\name{B}\times\name{A}=\name{i}$. For every $n\in\omega$ there is a maximal antichain $I_n$ which forces a truth value to the assertion $n\in\name{A}$. Set $U=\bigcup\{{\rm Dom}(q):q\in I_n,n\in\omega\}$. Since $|U|=\aleph_0$ we know that $\Vdash_{\mathbb{Q}}\name{B}\nsubseteq U$.

Consequently, there is some $\alpha<\mu$ such that $p_0\nVdash\alpha\notin\name{B}$, so one can choose a condition $p_1\geq p_0$ such that $p_1\Vdash\alpha\in\name{B}$. Without loss of generality, $\alpha\in\Dom(p_1)$. Let $p_2$ be $p_1\upharpoonright U$. Choose a natural number $n_*$ such that ${\rm sup}(\Dom(p_1(\alpha)))<n_*$. There are $n, p_3$ such that $p_2\leq p_3, \Dom(p_3)\subseteq U$ and $p_3\Vdash_{\mathbb{Q}}n_*<n\in\name{A}$. Define $p_4=p_1\upharpoonright(\Dom(p_1)\setminus U)\cup p_3$ and $p_5=p_4\cup\langle\alpha,1-\name{i}\rangle$. The condition $p_5$ forces $\name{c}(\alpha,n)\neq\name{i}$, a contradiction.

\hfill \qedref{ccohen}

So our problems are connected. We would like to show that under the continuum hypothesis there is a strong splitting family. We quote the following result (in a more general form), about strong polarized relations under the local assumption of the GCH. The proof appears in \cite{williams}:

\begin{proposition}
\label{ggch}
Polarized relation and the GCH. \newline 
Assume $2^\kappa = \kappa^+$. \newline 
\Then\ $\binom{\kappa^+}{\kappa} \nrightarrow \binom{\kappa^+}{\kappa}^{1,1}_2$.
\end{proposition}

\par \noindent We can conclude:

\begin{corollary}
\label{ffdirection}
The existence of strong splitting families. \newline 
Suppose $2^{\aleph_0} = \aleph_1$. \newline 
\Then\ there exists a strong splitting family.
\end{corollary}

\par \noindent \emph{Proof}. \newline 
By proposition \ref{ggch}, $\binom{\aleph_1}{\aleph_0} \nrightarrow \binom{\aleph_1}{\aleph_0}^{1,1}_2$. Since $\aleph_1 \leq \mathfrak{s} \leq 2^{\aleph_0}$, we have (under the continuum hypothesis) $\mathfrak{s} = \aleph_1$, so $\binom{\mathfrak{s}}{\omega} \nrightarrow \binom{\mathfrak{s}}{\omega}^{1,1}_2$, and by claim \ref{eequiv} we know that a strong splitting family exists.

\hfill \qedref{ffdirection}

\begin{definition}
\label{ssn}
The strong splitting numbers.
\begin{enumerate}
\item [$(a)$] Let $\mathfrak{ss}$ (= the strong splitting number) be the minimal cardinality of a strong splitting family, if one exists.
\item [$(b)$] Similarly, $\mathfrak{ss'}$ (= the very strong splitting number) is the minimal cardinality of a very strong splitting family, if one exists.
\end{enumerate}
\end{definition}

\begin{claim}
\label{sssequal}
The equality claim. \newline 
If there is a strong splitting family, then $\mathfrak{s} = \mathfrak{ss}$.
\end{claim}

\par \noindent \emph{Proof}. \newline 
Let $\kappa$ be $\mathfrak{ss}$. By remark \ref{ttrrivial}, $\kappa \geq \mathfrak{s}$. Assume toward contradiction that $\kappa > \mathfrak{s}$. We shall prove that $\binom{\mathfrak{s}}{\omega} \rightarrow \binom{\mathfrak{s}}{\omega}^{1,1}_2$. Let $c : \mathfrak{s} \times \omega \rightarrow 2$ be any coloring. For every $\alpha < \mathfrak{s}$, define $S_\alpha = \{ n \in \omega : c(\alpha, n) = 0 \}$. Define $\mathcal{F}_c \equiv \mathcal{F} = \{ S_\alpha : \alpha < \mathfrak{s} \}$. Since $|\mathcal{F}| \leq \mathfrak{s} < \kappa$ (by the assumption toward contradiction), $\mathcal{F}$ is not a strong splitting family. 

Choose a witness, $B$. It means that $B \in [\omega]^\omega$ and $|\mathcal{F}_B| = \mathfrak{s}$. For each $S_\alpha \in \mathcal{F}_B$ we have 
$(B \subseteq^* S_\alpha)$ or $(B \subseteq^* \omega \setminus S_\alpha)$. Without loss of generality, $B \subseteq^* S_\alpha$ for every $\alpha < \mathfrak{s}$, and moreover, $B \subseteq S_\alpha$ for $\mathfrak{s}$-many $\alpha$'s (recall that $\cf(\mathfrak{s}) > \aleph_0$). Let $H$ be the set $\{ \alpha < \mathfrak{s} : S_\alpha \in \mathcal{F}_B \}$. By the construction, $c \upharpoonright (H \times B) \equiv 0$, so $\binom{\mathfrak{s}}{\omega} \rightarrow \binom{\mathfrak{s}}{\omega}^{1,1}_2$.
By Claim \ref{eequiv}, $\binom{\mathfrak{s}}{\omega} \rightarrow \binom{\mathfrak{s}}{\omega}^{1,1}_2$ implies that there is no strong splitting family, a contradiction.

\hfill \qedref{sssequal}

\begin{remark}
\label{milovic}
We thank David Milovich for informing us the consistency of $\mathfrak{ss}=\kappa$ for every regular $\kappa>\aleph_0$. The interested reader can find the proof (among other results) in \cite{MR2467217}.
\end{remark}

It is not known if $\mathfrak{s}$ can be a singular cardinal. In general, it seems that $\mathfrak{ss}$ is more convenient to deal with. By the previous results, under some circumstances, we can infer about $\mathfrak{s}$ from what we know about $\mathfrak{ss}$. The following claim illustrates this idea:

\begin{claim}
\label{rregul}
Suppose there is a strong splitting family. \newline 
If $\mathfrak{s} = \mathfrak{c}$, \then\ $\mathfrak{s}$ is a regular cardinal. 
Similarly, if a very strong splitting family exists then $\mathfrak{ss'}=\mathfrak{c}$ implies the regularity of $\mathfrak{c}$.
\end{claim}

\par \noindent \emph{Proof}. \newline 
Assume toward contradiction that $\cf(\mathfrak{s}) < \mathfrak{s}$. By the assumption of the claim, $\mathfrak{c}$ is also a singular cardinal. Choose an unbounded increasing sequence of ordinals $\langle \xi_\gamma : \gamma < \cf(\mathfrak{c}) \rangle$, whose limit is $\mathfrak{c}$. Choose a strong splitting family $\mathcal{F} = \{ S_\alpha : \alpha < \mathfrak{s} \}$.

For every $B \in [\omega]^\omega$ we know that $|\mathcal{F}_B| < \mathfrak{s}$. Denote the set $\{ \alpha < \mathfrak{s} : S_\alpha \in \mathcal{F}_B \}$ by $H_B$. Clearly, $|H_B| < \mathfrak{s}$. We claim that there is an $H \subseteq \mathfrak{c}, |H|=\cf(\mathfrak{c})$ such that $H \nsubseteq H_B$ for every $B \in [\omega]^\omega$.

Choose an enumeration $\{ H_\alpha : \alpha < \mathfrak{c} \}$ of the $H_B$'s. We can assume that ${\rm sup} \{ |H_\varepsilon| : \varepsilon < \xi_\gamma \} < \mathfrak{c}$ for every $\gamma < \cf(\mathfrak{c})$. For every $\gamma < \cf(\mathfrak{c})$ choose $a_\gamma \in \mathfrak{c} \setminus (\bigcup \{ H_\varepsilon : \varepsilon < \xi_\gamma \} \cup \{ a_\beta : \beta < \gamma \})$. This is possible, since $|\bigcup \{ H_\varepsilon : \varepsilon < \xi_\gamma \}| \leq |\xi_\gamma| \cdot {\rm sup} \{ |H_\varepsilon| : \varepsilon < \xi_\gamma \} < \mathfrak{c}$. Set $H = \{ a_\gamma : \gamma < \cf(\mathfrak{c}) \}$, and we have the desired $H$.

Now define $\mathcal{F}' = \{ S_\alpha : \alpha \in H \}$. Clearly, $|\mathcal{F}'| = |H| \leq \cf(\mathfrak{c}) < \mathfrak{c} = \mathfrak{s}$. But if $B \in [\omega]^\omega$ then there is an ordinal $\alpha \in H$ so that $S_\alpha$ splits $B$ (by the fact that $H \nsubseteq H_B$). So $\mathcal{F}'$ is a splitting family whose cardinality is strictly less than $\mathfrak{s}$, a contradiction. The proof of the second assertion is identical.

\hfill \qedref{rregul}

\medskip 

Recall that ${\rm cov}(\lambda, \mu, \theta, 2)$ is the minimal cardinality of a family of subsets of $\lambda$, the cardinality of each member is below $\mu$, such that every set in $[\lambda]^{< \theta}$ is covered by a member from this family. By a similar argument, we can conclude:

\begin{corollary}
\label{singgg} Splitting properties and covering numbers.
\begin{enumerate}
\item [$(a)$] Suppose ${\rm cov}(\mathfrak{s}, \mathfrak{s}, \mathfrak{s}, 2) > \mathfrak{c}$. \newline 
\Then\ there is no strong splitting family.
\item [$(b)$] If $\mu>\cf(\mu)$ and ${\rm cov}(\mu,\mu,\mu,2) > \mathfrak{c}$ then there is no very strong splitting family of cardinality $\mu$.
\end{enumerate}
\end{corollary}

\par \noindent \emph{Proof}. \newline 
Assume toward contradiction that there exists a strong splitting family, and let $\mathcal{F} = \{ S_\alpha : \alpha < \mathfrak{s} \}$ exemplify it. For every $B \in [\omega]^\omega$ let $H_B = \{ \alpha < \mathfrak{s} : S_\alpha \in \mathcal{F}_B \}$. By our assumption toward contradiction we know that $|H_B| < \mathfrak{s}$ for every $B$.

By the assumption that ${\rm cov}(\mathfrak{s}, \mathfrak{s}, \mathfrak{s}, 2) > \mathfrak{c}$, one can pick a set $H \subseteq \mathfrak{s}, |H| < \mathfrak{s}$ which is not covered by the family $\{ H_B : B \in [\omega]^\omega \}$. Now set $\mathcal{F}' = \{ S_\alpha : \alpha \in H \}$.

Since $|H| < \mathfrak{s}$, we know that $|\mathcal{F}'| < \mathfrak{s}$. By the nature of $H$ we have $H \nsubseteq H_B$ for every $B \in [\omega]^\omega$. Hence, given $B \in [\omega]^\omega$ one can pick an ordinal $\alpha \in H \setminus H_B$, so $S_\alpha$ splits $B$. It means that $\mathcal{F}'$ is a splitting family whose cardinality is strictly less than $\mathfrak{s}$, a contradiction. Here, again, the second part of the corollary follows in a similar way.

\hfill \qedref{singgg}

\begin{remark}
\label{hhighs}
The definition of $\mathfrak{s}$ is generalized naturally to higher cardinals (see, for example, \cite{MR1450512}). $\mathfrak{s}_\lambda$ is the minimal cardinality of a $\lambda$-splitting family in $[\lambda]^\lambda$. It is known that $\mathfrak{s}_\lambda > \lambda$ iff $\lambda$ is weakly compact (see \cite{MR1450512}). 

In this case, the main claim of this paper can be applied to $\mathfrak{s}_\lambda$, yielding $\binom{\mathfrak{s}_\lambda}{\lambda} \nrightarrow \binom{\mathfrak{s}_\lambda}{\lambda}^{1,1}_2$ iff there is a strong $\lambda$-splitting family in $[\lambda]^\lambda$. The existence result under the assumption $2^\lambda = \lambda^+$ follows. 
\end{remark}

\par \noindent In a subsequent paper (see \cite{GaSh995}) we prove that the positive relation $\binom{\mathfrak{s}}{\omega} \rightarrow \binom{\mathfrak{s}}{\omega}^{1,1}_2$ (and hence the non-existence of strong splitting families) is also consistent with ZFC.

\newpage

\section{Polarized relations and the weak diamond}
We prove, in this section, the consistency of positive polarized relations with the weak diamond. The following definition comes from \cite{MR0469756}:

\begin{definition}
\label{weakdiamond} The weak diamond. \newline
$\Phi_{\aleph_1}$ means that for every $F:{^{\omega_1>}2} \rightarrow 2$ there exists $g\in {^{\omega_1}2}$ so that for every $f\in {^{\omega_1}2}$ the set $S=\{\alpha<\aleph_1 : F(f\upharpoonright\alpha)=g(\alpha)\}$ is stationary (in $\aleph_1$).
\end{definition}

The idea is pretty simple. The diamond sequence provides a tool for guessing many initial segments of every $A\subseteq\aleph_1$ (in the sense of $A\cap\alpha$ for stationarily many $\alpha$-s). The weak diamond does not give the set $A$, but it gives a way to guess the color of $A\cap\alpha$ (again, for stationarily many $\alpha$-s) once a coloring of ${^{\omega_1>}2}$ is in hand.

It is shown in \cite{MR0469756} that $2^{\aleph_0} < 2^{\aleph_1} \Rightarrow \Phi_{\aleph_1}$ (and as noted by Uri Abraham, $\Phi_{\aleph_1} \Rightarrow 2^{\aleph_0} < 2^{\aleph_1}$, so actually we have an equivalence). This gives rise to the following simple fact:

\begin{proposition}
\label{wweakcf} Weak diamond and low cofinality. \newline
Suppose $\cf(2^\omega)=\aleph_1$. $\then\ \Phi_{\aleph_1}$ holds.
\end{proposition}

\par \noindent \emph{Proof}. \newline 
By Zermelo-K\"onig, $\cf(2^{\omega_1})>\aleph_1$, hence $\cf(2^\omega)=\aleph_1$ implies $2^\omega\neq 2^{\omega_1}$, i.e., $2^\omega<2^{\omega_1}$ which yields $\Phi_{\aleph_1}$.

\hfill \qedref{wweakcf}

As described in the introduction, one could suspect that $\Phi_{\aleph_1}$ entails negative polarized relations (similar to the impact of the real diamond). The following claims show that $\Phi_{\aleph_1}$ is not strong enough. The first claim deals with $\binom{\omega_1}{\omega} \rightarrow \binom{\omega_1}{\omega}^{1,1}_n$, and the second deals with $\binom{\mathfrak{c}}{\omega} \rightarrow \binom{\mathfrak{c}}{\omega}^{1,1}_n$.

We shall use the Mathias forcing $\mathbb{M}_D$ for proving the main result of this section.
Let $D$ be a nonprincipal ultrafilter on $\omega$. We define $\mathbb{M}_D$ as follows. The conditions in $\mathbb{M}_D$ are pairs of the form $(s, A)$ when $s \in [\omega]^{< \omega}$, $A \in D$ and ${\rm min}(A)>{\rm max}(s)$. For the order, $(s_1, A_1) \leq (s_2, A_2)$ iff $s_1 \subseteq s_2, A_1 \supseteq A_2$ and $s_2 \setminus s_1 \subseteq A_1$. The Mathias forcing is $\sigma$-centered, hence satisfies the $ccc$. It follows that a finite support iteration of these forcing notions is also $ccc$.

Let $G \subseteq \mathbb{M}_D$ be generic over ${\rm \bf V}$. The Mathias real $x_G$ is defined as $\bigcup \{ s : \exists A \in D, (s, A) \in G \}$. Notice that in ${\rm \bf V}^{\mathbb{M}_D}$ we have $(x_G \subseteq^* B) \vee (x_G \subseteq^* \omega \setminus B)$ for every $B \in [\omega]^\omega\cap {\rm \bf V}$. In a way, it means that the Mathias forcing adds a pseudo-intersection to the ultrafilter $D$. By iterating these forcing notions we create the monochromatic subsets.

\begin{claim}
\label{wwwandwww} Weak diamond and $\omega_1$. \newline 
The positive relation $\binom{\omega_1}{\omega} \rightarrow \binom{\omega_1}{\omega}^{1,1}_n$ is consistent with $\Phi_{\aleph_1}$, and even with $\cf(2^{\aleph_0})=\aleph_1$.
\end{claim}

\par \noindent \emph{Proof}. \newline 
By theorem \ref{mmt} below.

\hfill \qedref{wwwandwww}

\begin{theorem}
\label{mmt} Weak diamond and the continuum.
\begin{enumerate}
\item [$(a)$] The strong relation $\binom{\mathfrak{c}}{\omega} \rightarrow \binom{\mathfrak{c}}{\omega}^{1,1}_n$ is consistent with $\Phi_{\aleph_1}$, and even with $\cf(2^{\aleph_0})=\aleph_1$.
\item [$(b)$] Moreover, suppose $\aleph_1\leq\theta=\cf(\theta)\leq\mu$, and $\mu=\mu^{\aleph_0}$ (in ${\rm \bf V}$). There is a $ccc$ forcing notion $\mathbb{P}, |\mathbb{P}|=\mu, \Vdash_{\mathbb{P}}2^{\aleph_0}=\mu$, and for every $\lambda\in(\aleph_0,\mu]$ if $\cf(\lambda)\notin \{\aleph_0,\theta\}$ then ${\rm \bf V}^\mathbb{P} \models \binom{\lambda}{\omega} \rightarrow \binom{\lambda}{\omega}^{1,1}_n$.
\end{enumerate}
\end{theorem}

\par \noindent \emph{Proof}. \newline 
We prove the second assertion of the theorem. By choosing any $\mu=\mu^{\aleph_0}>\cf(\mu)=\aleph_1$ and $\theta=\aleph_2, \lambda=\mu$ we will get a proof to the first assertion.
So choose an ordinal $\delta\in[\mu,\mu^+)$ so that $\cf(\delta)=\theta$. We define a finite support iteration $\langle \mathbb{P}_i, \name{\mathbb{Q}}_j : i \leq \delta, j < \delta \rangle$ of $ccc$ forcing notions, such that $|\mathbb{P}_i| = \mu$ for every $i \leq \delta$.

Let $\name{\mathbb{Q}}_0$ be (a name of) a forcing notion which adds $\mu$ reals (e.g., Cohen forcing). For every $j < \omega_1$ let $\name{D}_j$ be a $\mathbb{P}_j$-name of a nonprincipal ultrafilter on $\omega$ over the extension with $\mathbb{P}_j$. Let $\name{\mathbb{Q}}_{1+j}$ be the Mathias forcing $\mathbb{M}_{\name{D}_j}$. $\name{\mathbb{Q}}_{1+j}$ is a $ccc$ forcing notion which adds an infinite set $\name{A}_j \subseteq \omega$ such that $(\forall B \in \name{D}_j) (\name{A}_j \subseteq^* B \vee \name{A}_j \subseteq^* \omega \setminus B)$. At the end, set $\mathbb{P} = \bigcup \{\mathbb{P}_i : i < \delta\}$.

Since every component satisfies the $ccc$, and we use finite support iteration, $\mathbb{P}$ is also a $ccc$ forcing notion and hence no cardinal is collapsed in ${\rm \bf V}^\mathbb{P}$. By virtue of $\name{\mathbb{Q}}_0$, $2^{\aleph_0} = \mu$ after forcing with $\mathbb{P}$ (since $\name{\mathbb{Q}}_0$ adds $\mu$-many reals, and the length of the iteration is $\delta$, which is of size $\mu$). Let $\lambda$ be any cardinal in $(\aleph_0,\mu]$ such that $\cf(\lambda)\notin\{\aleph_0,\theta\}$, and let $n$ be a finite ordinal. Our goal is to prove that $\binom{\lambda}{\omega} \rightarrow \binom{\lambda}{\omega}^{1,1}_n$ in ${\rm \bf V}^\mathbb{P}$.

Let $\name{c}$ be a name of a function from $\lambda \times \omega$ into $n$. For every $\alpha < \lambda$ we have a name (in ${\rm \bf V}$) to the restriction $\name{c} \upharpoonright (\{ \alpha \} \times \omega)$. $\mathbb{P}$ is $ccc$, hence the color of every pair of the form $(\alpha, n)$ is determined by an antichain which includes at most $\aleph_0$ conditions. Since we have to decide the color of $\aleph_0$-many pairs in $\name{c} \upharpoonright (\{ \alpha \} \times \omega)$, and the length of $\mathbb{P}$ is $\delta$, $\cf(\delta)\geq\aleph_1$ we know that $\name{c} \upharpoonright (\{ \alpha \} \times \omega)$ is a name in $\mathbb{P}_{i(\alpha)}$ for some $i(\alpha) < \delta$.

For every $j < \delta$ let $\mathcal{U}_j$ be the set $\{ \alpha < \lambda : i(\alpha) \leq j \}$, so $\langle\mathcal{U}_j:j<\delta\rangle$ is $\subseteq$-increasing with union $\lambda$. Recall that $\cf(\delta)=\theta \neq \cf(\lambda)$, hence for some $j < \delta$ we have $\mathcal{U}_j \in [\lambda]^\lambda$. Choose such $j$, and denote $\mathcal{U}_j$ by $\mathcal{U}$. We shall try to show that $\mathcal{U}$ can serve (after some shrinking) as the first coordinate in the monochromatic subset. 

Choose a generic subset $G \subseteq \mathbb{P}$, and denote $\name{A}_j[G]$ by $A$. For each $\alpha \in \mathcal{U}$ we know that $\name{c} \upharpoonright (\{ \alpha \} \times A)$ is constant, except a possible mistake over a finite subset of $A$. But this mistake can be amended.

For every $\alpha \in \mathcal{U}$ choose $k(\alpha) \in \omega$ and $m(\alpha) < n$ so that $(\forall \ell \in A)[\ell \geq k(\alpha) \Rightarrow \name{c}[G](\alpha, \ell) = m(\alpha)]$. $n$ is finite and $\cf(\lambda) > \aleph_0$, so one can fix some $k \in \omega$ and a color $m < n$ such that for some $\mathcal{U}_1 \in [\mathcal{U}]^\lambda$ we have $\alpha \in \mathcal{U}_1 \Rightarrow k(\alpha) = k \wedge m(\alpha) = m$.

Let $B$ be $A \setminus k$, so $B \in [\omega]^\omega$. By the fact that $\mathcal{U}_1 \subseteq \mathcal{U}$ we know that $c(\alpha, \ell) = m$ for every $\alpha \in \mathcal{U}_1$ and $\ell \in B$, so $\mathcal{U}_1\times B$ is monochromatic under $c$, yielding the positive relation $\binom{\lambda}{\omega} \rightarrow \binom{\lambda}{\omega}^{1,1}_n$, as required.

\hfill \qedref{mmt}

\begin{remark}
\label{sstt}
Assume $\lambda$ is an uncountable regular cardinal. Denote by $\binom{\lambda}{\mu} \rightarrow_{\rm st} \binom{\lambda}{\mu}^{1,1}_n$ the assertion that for every coloring $c : \lambda \times \mu \rightarrow n$ there exists $B \in [\mu]^\mu$ and a stationary subset $\mathcal{U} \subseteq \lambda$ so that $c \upharpoonright \mathcal{U} \times B$ is constant. Our proof gives the consistency of $\binom{\lambda}{\omega} \rightarrow_{\rm st} \binom{\lambda}{\omega}^{1,1}_n$, when $\lambda$ is regular.
\end{remark}

Recall that if $\cf(\kappa)>\aleph_0$ and $\kappa<\mathfrak{s}$ then $\binom{\kappa}{\omega}\rightarrow\binom{\kappa}{\omega}^{1,1}_2$ (this is claim 1.4 in \cite{gash964}). We can use this claim for accomplishing the proof of \ref{wwwandwww}, as shown below:

\par \noindent \emph{Proof of \ref{wwwandwww}}: \newline 
Choose $\theta=\aleph_2\leq\mu, \cf(\mu)=\aleph_1$ and $\mu=\mu^{\aleph_0}$. Use the iteration in the proof of \ref{mmt} over some ordinal $\delta\in(\mu,\mu^+)$ so that $\cf(\delta)=\aleph_2$. By the properties of the Mathias forcing we have $\mathfrak{s}=\theta=\aleph_2$ in the extension. For showing this, we shall prove that $\aleph_1<\mathfrak{s}\leq\aleph_2$.

If $\{S_\alpha:\alpha<\kappa\}$ is a splitting family in ${\rm \bf V}^\mathbb{P}$ and $\kappa<\aleph_2$, then there exists an ordinal $j<\delta$ such that $S_\alpha\in {\rm \bf V}^{\mathbb{P}_j}$ for every $\alpha<\kappa$ (since $\cf(\delta)=\omega_2$). But then, the Mathias real added in the $j$-th stage is almost included in $S_\alpha$ or its complement for every $\alpha<\kappa$, a contratiction. Hence $\aleph_1<\mathfrak{s}$.

On the other hand, one can introduce a splitting family of size $\aleph_2$ in ${\rm \bf V}^\mathbb{P}$. Along a finite support iteration, a Cohen real is added at every limit stage. Choose a cofinal sequence of limit ordinals in $\delta$, of length $\aleph_2$. Since every Cohen real is a splitting real (over the old universe), the collection of the Cohen reals along the cofinal sequence establishes a splitting family of size $\aleph_2$, so $\mathfrak{s}\leq\aleph_2$ in ${\rm \bf V}^\mathbb{P}$.

It follows from the remark above that $\binom{\omega_1}{\omega} \rightarrow \binom{\omega_1}{\omega}^{1,1}_2$ since $\omega_1=\cf(\omega_1)<\mathfrak{s}$. On the other hand, $\mu=2^{\aleph_0}$ so $\cf(2^{\aleph_0})=\aleph_1$, hence $\Phi_{\aleph_1}$ as required.

\hfill \qedref{wwwandwww}

\newpage

\bibliographystyle{amsplain}
\bibliography{arlist}

\end{document}